\documentclass[a4paper]{amsart}

\usepackage{amssymb}

\newtheorem{proposition}[equation]{Proposition}

\numberwithin{equation}{section}

\theoremstyle{definition}

\newtheorem*{example*}{Example}

\newtheorem*{remark*}{Remark}

\newcommand{\bF}{{\mathbb F}}
\newcommand{\bZ}{{\mathbb Z}}

\newcommand{\cV}{{\mathcal V}}
\newcommand{\frg}{{\mathfrak g}}

\newcommand{\frf}{{\mathfrak f}}
\newcommand{\fre}{{\mathfrak e}}

\newcommand{\calT}{{\mathcal T}}

\newcommand{\calL}{{\mathcal L}}

\newcommand{\calO}{{\mathcal O}}

\newcommand{\subo}{_{\bar 0}}
\newcommand{\subuno}{_{\bar 1}}

\newcommand{\bil}{{\textup{b}}}

\DeclareMathOperator{\eespan}{span}
\providecommand{\espan}[1]{\eespan\left\{ #1\right\}}

 \newcommand{\tri}{\mathfrak{tri}}

 \newcommand{\frsl}{{\mathfrak{sl}}}
 \newcommand{\frsp}{{\mathfrak{sp}}}
 \newcommand{\frso}{{\mathfrak{so}}}
 \newcommand{\frpsl}{{\mathfrak{psl}}}
 \newcommand{\frgl}{{\mathfrak{gl}}}
 \newcommand{\frpgl}{{\mathfrak{pgl}}}
 \newcommand{\frosp}{{\mathfrak{osp}}}

  \newcommand{\frbr}{{\mathfrak{br}}}
   \newcommand{\frel}{{\mathfrak{el}}}

 \DeclareMathOperator{\tr}{tr}
  
 \DeclareMathOperator{\ad}{ad}
 
 \DeclareMathOperator{\der}{\mathfrak{der}}
 
 \DeclareMathOperator{\inder}{\mathfrak{inder}}
 
 \DeclareMathOperator{\End}{End}
 
 \DeclareMathOperator{\Mat}{Mat}

 \DeclareMathOperator{\charac}{char}

\def\bigstrut{\vrule height 14pt width 0ptdepth 2pt}

\begin{document}

\title{Tits construction of the exceptional simple Lie algebras}

\author[Alberto Elduque]{Alberto Elduque$^{\star}$}
 \thanks{$^{\star}$ Supported by the Spanish Ministerio de
 Educaci\'on y Ciencia
 and FEDER (MTM 2007-67884-C04-02) and by the
Diputaci\'on General de Arag\'on (Grupo de Investigaci\'on de
\'Algebra)}
 \address{Departamento de Matem\'aticas e
 Instituto Universitario de Matem\'aticas y Aplicaciones,
 Universidad de Zaragoza, 50009 Zaragoza, Spain}
 \email{elduque@unizar.es}

\dedicatory{Dedicated to Professor Jacques Tits on the occasion of his eightieth birthday}

\date{\today}

\subjclass[2000]{Primary: 17B60. Secondary: 17B50, 17B25, 17C50, 17A75}

\keywords{Tits construction, Lie algebra, exceptional, simple, Jordan, composition}

\begin{abstract}
The classical Tits construction of the exceptional simple Lie algebras has been extended in a couple of directions by using either Jordan superalgebras or composition superalgebras. These extensions are reviewed here. The outcome has been the discovery of some new simple modular Lie superalgebras.
\end{abstract}

\maketitle


\section{Introduction}\label{se:Introduction}

In 1966, Tits provided a beautiful unified construction of the exceptional simple Lie algebras $F_4$, $E_6$, $E_7$ and $E_8$ over fields of characteristic $\ne 2,3$ \cite{Tits66}. This construction depends on two algebras: an alternative algebra of degree $2$ and a Jordan algebra of degree $3$. The most interesting cases appear when semisimple alternative algebras and central simple Jordan algebras are considered in this construction. Then the alternative algebra becomes a composition algebra, and hence its dimension is restricted to $1$, $2$, $4$ or $8$, while the Jordan algebra becomes a form of the Jordan algebra of hermitian $3\times 3$ matrices over a second composition algebra.  Freudenthal's Magic Square (\cite{Freu64}) is obtained with these ingredients.

\smallskip

In recent years, some extensions of Tits construction have been considered. Benkart and Zelmanov \cite{BZ96} considered Lie algebras graded by root systems. They realized that the Lie algebras graded over either the root system $G_2$ or $F_4$ could be described in terms of a generalized Tits construction, and that some specific simple Jordan superalgebras (of dimension $3$ and $4$) could be plugged into this Tits construction, instead of Jordan algebras (see also \cite{BE03}). The outcome is that Freudenthal's Magic Square can be enlarged to a rectangle that includes, in characteristic $\ne 2,3,5$, the exceptional simple Lie superalgebras $D(2,1;t)$, $G(3)$ and $F(4)$ in Kac's classification \cite{Kac77}. Moreover, in characteristic $5$ another simple Jordan superalgebra: Kac's ten-dimensional algebra, can be used too, and this lead to the discovery of a new simple modular Lie superalgebra specific of this characteristic \cite{EldModular}´, whose even part is of type $B_5$ and the odd part is the spin module for the even part.

\smallskip

In a different direction, Freudenthal's Magic Square presents a symmetry which is not obvious from Tits construction. More symmetric constructions have been considered by different authors (\cite{V}, \cite{AF}, \cite{LM1}, \cite{LM2}, \cite{BS1}, \cite{BS2} or \cite{EldIbero1}). These constructions are based on two composition algebras and their triality Lie algebras, and they are symmetric on the two composition algebras involved, whence the symmetry of the outcome: the Magic Square. Besides, these construction are valid too over fields of characteristic $3$.

But over fields of characteristic $3$ there are composition superalgebras with nontrivial odd part, as this phenomenon is specific of characteristics $2$ and $3$ (in characteristic $2$, these superalgebras are just composition algebras with a grading over $\bZ_2$). These composition superalgebras in characteristic $3$ can be plugged into the symmetric construction of the Magic Square. The outcome is a larger square, which extends the Magic Square with the addition of two rows and columns filled with (mostly) simple Lie superalgebras. With one exception, the Lie superalgebras that appear have no counterpart in Kac's classification. There are $10$ such new simple modular Lie superalgebras in characteristic $3$.

\smallskip

Moreover, it turns out that the simple Lie superalgebra in characteristic $5$ mentioned above, obtained by means of the Tits construction with ingredients a Cayley algebra and the ten-dimensional Kac Jordan superalgebra, and the $10$ simple Lie superalgebras that appear in the larger square in characteristic $3$ almost exhaust (there are just $3$ other superalgebras) the list of exceptional simple modular Lie superalgebras with a Cartan matrix in characteristic $>2$ in the recent classification by Bouarroudj, Grozman and Leites \cite{BGLgordo}.

\medskip

The aim of this paper is to review these two directions in which the Tits construction has been extended.

In the next section, the classical Tits construction will be recalled. Then a section will be devoted to see how some Jordan superalgebras can be used instead of Jordan algebras. This will result in a larger \emph{Supermagic Rectangle} which contains the exceptional simple classical Lie superalgebras in Kac's classification and a new simple Lie superalgebra in characteristic $5$.

The fourth section will be devoted to explain the more symmetric construction of Freudenthal's Magic Square in terms of two composition algebras and then the fifth section will show how this symmetric construction allows, in characteristic $3$, to enlarge the Magic Square to a \emph{Supermagic Square} with a whole family of new simple Lie superalgebras which appear only in this characteristic. Some related comments and remarks will be given in a final section.

\smallskip

No proofs will be given, but references will be provided throughout. For the sake of simplicity, all the vector spaces and algebras and superalgebras considered will be defined over an algebraically closed ground field $\bF$ of characteristic $\ne 2$.

\bigskip


\section{Tits construction}\label{se:Tits}

As mentioned in the introduction, in 1966 Tits (\cite{Tits66}) gave a unified construction of the exceptional simple classical Lie algebras by means of two
ingredients: a composition algebra and a degree three simple
Jordan algebra. Here this construction will be reviewed following closely
\cite[\S 1]{EldTits3} and the approach used by Benkart and Zelmanov in
\cite{BZ96}.

Let $C$ be a composition algebra over the ground field $\bF$
with norm $n$. Thus, $C$ is a finite dimensional unital algebra over $\bF$,
endowed with the nondegenerate quadratic form $n:C\rightarrow \bF$ such that
$n(ab)=n(a)n(b)$ for any $a,b\in C$. Then, each element satisfies
the degree $2$ equation
\[
a^2-t_C(a)a+n(a)1=0,
\]
where $t_C(a)=n(a,1)\,\bigl(=n(a+1)-n(a)-1\bigr)$ is called the
\emph{trace}. The subspace of trace zero elements will be denoted by
$C^0$.

Moreover, for any $a,b\in C$, the linear map $D_{a,b}:C\rightarrow
C$ given by
\[
D_{a,b}(c)=[[a,b],c]+3(a,c,b)
\]
where $[a,b]=ab-ba$ is the commutator, and $(a,c,b)=(ac)b-a(cb)$ the
associator, is a derivation: the \emph{inner derivation} determined
by the elements $a,b$ (see \cite[Chapter III, \S 8]{Schafer}). These
derivations satisfy
\[
D_{a,b}=-D_{b,a},\quad D_{ab,c}+D_{bc,a}+D_{ca,b}=0,
\]
for any $a,b,c\in C$. The linear span of these derivations will be
denoted by $\inder C$. It is an ideal of the whole Lie algebra of
derivations $\der C$ and, if the characteristic is $\ne 3$, it is
the whole $\der C$.

The dimension of $C$ is restricted to $1$, $2$, $4$ (quaternion
algebras) and $8$ (Cayley algebras). Over our
algebraically closed field $\bF$, the only composition algebras are, up to
isomorphism, the ground field $\bF$, the cartesian product of two
copies of the ground field $\bF\times \bF$, the algebra of two by two
matrices $\Mat_2(\bF)$ (where the norm of a matrix is its determinant), and the split Cayley algebra $C(\bF)$. (See, for
instance, \cite[Chapter 2]{ZSSS}.)

\smallskip

Now, let $J$ be a unital Jordan algebra with a \emph{normalized
trace} $t_J:J\rightarrow \bF$. That is, $J$ is a commutative algebra over $\bF$ which satisfies the Jordan identity:
\[
x^2(yx)=(x^2y)x
\]
for any $x,y\in J$, and $t_J$ is a linear map such
that $t_J(1)=1$ and $t_J\bigl((xy)z\bigr)=t_J\bigl(x(yz)\bigr)$ for
any $x,y,z\in J$. The archetypical example of a Jordan algebra is the subspace of symmetric elements in an associative algebra with an involution. This subspace is not closed under the associative product, but it is indeed closed under the symmetrized product $\frac{1}{2}(xy+yx)$. With this symmetrized product, it becomes a Jordan algebra. (One may consult \cite{Jac68} for the main properties of finite dimensional Jordan algebras.)

Then, for such a unital Jordan algebra, there is the decomposition $J=\bF 1\oplus J^0$, where $J^0=\{x\in J:
t_J(x)=0\}$. For any $x,y\in J^0$, the product $xy$ splits as
\[
xy=t_J(xy)1+x*y,
\]
with $x*y\in J^0$. Then $x*y=xy-t_J(xy)1$ gives a commutative
multiplication on $J^0$. The linear map $d_{x,y}:J\rightarrow J$
defined by
\[
d_{x,y}(z)=x(yz)-y(xz),
\]
is the \emph{inner derivation} of $J$ determined by the elements $x$ and
$y$. Since $d_{1,x}=0$ for any $x$, it is enough to deal with the
inner derivations $d_{x,y}$, with $x,y\in J^0$. The linear span of
these derivations will be denoted by $\inder J$, which is an ideal
of the whole Lie algebra of derivations $\der J$.

\smallskip

Given $C$ and $J$ as before, consider the space
\[
\calT(C,J)=\inder C\oplus \bigl(C^0\otimes J^0\bigr)\oplus
\inder J
\]
(unadorned tensor products are always considered over $\bF$), with the
anticommutative multiplication $[.,.]$ specified by:
\begin{equation}\label{eq:TCJproduct}
\begin{split}
\bullet&\ \textrm{$\inder C$ and $\inder J$ are Lie subalgebras,}\\
\bullet&\ [\inder C,\inder J]=0,\\
\bullet&\ [D,a\otimes x]=D(a)\otimes x,\ [d,a\otimes x]=a\otimes
d(x),\\
\bullet&\ [a\otimes x,b\otimes y]=t_J(xy)D_{a,b}+\bigl([a,b]\otimes
x*y\bigr)+2t_C(ab)d_{x,y},
\end{split}
\end{equation}
for all $D\in \inder C$, $d\in \inder J$, $a,b\in C^0$, and
$x,y\in J^0$.

\smallskip

The conditions for $\calT(C,J)$ to be a Lie algebra are the
following:

\[
\begin{split}
\textrm{(i)}&\quad  \displaystyle{\sum_{\circlearrowleft}
 t_C\bigl([a_{1}, a_{2}] a_{3}\bigr)\,
d_{(x_1 * x_2), x_3}}=0,\\[6pt]
\textrm{(ii)}&\quad  \displaystyle{\sum_{\circlearrowleft}
 t_J\bigl( (x_1 * x_2) x_{3}\bigr)
\,D_{[a_1, a_2], a_3}}=0,\\[6pt]
\textrm{(iii)}&\quad \displaystyle{\sum_{\circlearrowleft}
 \Bigl(D_{a_1,a_2}(a_3) \otimes t_J\bigl(x_1
x_2\bigr) x_3} + [[a_1, a_2],a_3] \otimes (x_1 * x_2)* x_3\\[-6pt]
  &\null\hspace{2in} +2
t_C(a_1 a_2) a_3\otimes d_{x_1, x_2}(x_3)\Bigr)=0,
\end{split}
\]

\noindent
for any $a_1,a_2,a_{3} \in C^0$ and any $x_1,x_2,x_3 \in J^0$. The
notation ``$\displaystyle{\sum_\circlearrowleft}$'' indicates
summation over the cyclic permutation of the variables.

These conditions appear in \cite[Proposition 1.5]{BE03}, but there
they are stated in the more general setting of superalgebras, a
setting we will deal with later on. In particular, over fields of
characteristic $\ne 3$, these conditions are fulfilled if $J$ is a
separable Jordan algebra of degree three over $\bF$ and
$t_J=\frac{1}{3}T$, where $T$ denotes the generic trace of $J$ (see
for instance \cite{Jac68}).

If the characteristic of our algebraically closed field $\bF$ is $\ne 3$,
the degree $3$ simple Jordan algebras are, up to isomorphism, the
algebras of $3\times 3$ hermitian matrices over a unital composition
algebra: $H_3(C')$ (see \cite{Jac68}). By varying $C$ and $C'$,
$\calT(C,H_3(C'))$ is a classical Lie algebra, and
Freudenthal's Magic Square (Table \ref{ta:FMS}) is obtained.

\begin{table}[h!]
\begin{center}
\begin{tabular}{c|cccc}
 $\calT(C,J)$&\bigstrut$H_3(\bF)$&$H_3(\bF\times\bF)$
 &$H_3(\Mat_2(\bF))$&$H_3(C(\bF))$\\
 \cline{1-5}
 \bigstrut$\bF$&$A_1$&$A_2$&$C_3$&$F_4$\\
 \bigstrut$\bF\times \bF$&$A_2$&$A_2\oplus A_2$&$A_5$&$E_6$\\
 \bigstrut$\Mat_2(\bF)$&$C_3$&$A_5$&$D_6$&$E_7$\\
 \bigstrut$C(\bF)$&$F_4$&$E_6$&$E_7$&$E_8$
\end{tabular}
\bigskip\null
\end{center}
\caption{Freudenthal's Magic Square}\label{ta:FMS}
\end{table}

Let us have a look at the rows in this Tits construction of
Freudenthal's Magic Square.

\medskip

\noindent\textbf{First row:}\quad Here $C=\bF$, so $C^0=0$ and
$\inder C=0$. Thus, $\calT(C,J)$ is just $\inder J$. In
particular, $\calT(\bF,J)$ makes sense and is a Lie algebra for any
Jordan algebra $J$.

\medskip

\noindent\textbf{Second row:}\quad Here $C=\bF\times \bF$, so $C^0$
consists of the scalar multiples of $(1,-1)$, and thus $\calT(C,J)$
can be identified with $J^0\oplus \inder J$. The elements in $J^0$
multiply as $[x,y]=4d_{x,y}$ because
$t_C\bigl((1,-1)^2\bigr)=t_C\bigl((1,1)\bigr)=2$. Given any unital Jordan
algebra with a normalized trace $J$, its \emph{Lie multiplication algebra} $\calL(J)$ (see
\cite{Schafer}) is the Lie subalgebra of the general linear Lie
algebra $\frgl(J)$ generated by $l_J=\{l_x:x\in J\}$, where
$l_x:y\mapsto xy$ denotes the left multiplication by $x$. Then
$\calL(J)=l_J\oplus\inder J$. The map
\[
\begin{split}
\calT(C,J)&\rightarrow \calL(J)\\
(1,-1)\otimes x+d&\mapsto 2l_x+d,
\end{split}
\]
is a monomorphism. Its image is the Lie subalgebra $\calL^0(J)=l_{J^0}\oplus \inder J$. Again this shows that $\calT(\bF\times \bF,J)$ makes sense and is a
Lie algebra for any Jordan algebra with a normalized trace. Given
any separable Jordan algebra of degree $3$ in
characteristic $\ne 3$, $\calL^0(J)$ is precisely the derived
algebra $[\calL(J),\calL(J)]$. This latter Lie algebra makes sense
for any Jordan algebra over any field. (Recall that the
characteristic is assumed to be $\ne 2$ throughout.)

\medskip

\noindent\textbf{Third row:}\quad Here $C=\Mat_2(\bF)$. Under these circumstances, $C^0$ is the simple
three-dimensional Lie algebra of $2\times 2$ trace zero matrices $\frsl_2(\bF)$ under the commutator $[a,b]=ab-ba$.
Besides, for any $a,b\in C^0$, the inner derivation $D_{a,b}$ is
just $\ad_{[a,b]}$, since $C$ is associative. Hence, $\inder C$
can be identified with $C^0$, and $\calT(C,J)$ with
\[
C^0\oplus \bigl(C^0\otimes J^0\bigr)\oplus\inder
J\simeq\bigl(C^0\otimes (\bF1\oplus J^0)\bigr)\oplus\inder J\simeq
\bigl(C^0\otimes J\bigr)\oplus \inder J,
\]
and the Lie bracket \eqref{eq:TCJproduct} in $\calT(C,J)$ becomes
the bracket in $\bigl(\frsl_2(\bF)\otimes J\bigr)\oplus \inder J$ given by
\[
\begin{split}
\bullet&\ \textrm{$\inder J$ is a Lie subalgebra,}\\
\bullet&\ [d,a\otimes x]=a\otimes d(x),\\
\bullet&\ [a\otimes x,b\otimes y]=\bigl([a,b]\otimes
xy\bigr)+2t_C(ab)d_{x,y},
\end{split}
\]
for any $a,b\in \frsl_2(\bF)$, $x,y\in J$, and $d\in\inder J$, since
$t_J(xy)1+x*y=xy$ for any $x,y\in J$. This bracket makes sense for
any Jordan algebra (not necessarily endowed with a normalized
trace), it goes back to \cite{Tits62}, and the resulting Lie algebra is the well-known Tits-Kantor-Koecher Lie algebra of the Jordan algebra $J$ (\cite{Koecher}, \cite{Kantor}).

\medskip

\noindent\textbf{Fourth row:}\quad In the last row, $C$ is a Cayley
algebra over $\bF$. If the characteristic of $\bF$ is
$\ne 3$, the Lie algebra $\der C=\inder C$ is a simple Lie
algebra of type $G_2$ (dimension $14$, and note that this is no longer true in characteristic $3$ \cite{AMEN}), and $C^0$ is its unique
seven dimensional irreducible module. In particular, the Lie
algebra $\calT\bigl(C(\bF),J\bigr)$ is a Lie algebra graded over the
root system $G_2$. These $G_2$-graded Lie algebras contain a simple
subalgebra isomorphic to $\der C(\bF)$ such that, as modules for
this subalgebra, they are direct sums of copies of modules of three
types: adjoint, the irreducible seven dimensional module, and the
trivial one dimensional module. These Lie algebras have been
determined in \cite{BZ96} and the possible Jordan algebras involved
are essentially the degree $3$ Jordan algebras.

\bigskip

\section{Jordan superalgebras}\label{se:JordanSuperalgebras}

In Kac's classification of the simple finite dimensional Lie superalgebras over an algebraically closed field of characteristic $0$, there appear three exceptional situations (that is, the superalgebras do not belong to families where the dimension grows indefinitely).
A description of these
simple Lie superalgebras can be found in Kac's seminal paper \cite{Kac77}.  Over an algebraically closed field $\bF$
of characteristic $\ne 2,3$, these exceptional Lie superalgebras
may be characterized by the following properties (see
\cite[Proposition 2.1.1]{Kac77}, where the characteristic is assumed to be $0$, but this restriction is not necessary here):
\begin{enumerate}

\item There is a unique $40$-dimensional simple classical Lie
superalgebra $F(4)$ for which the even part $F(4)\subo$ is a
Lie algebra of type $B_3\oplus A_1$, and the
representation of $F(4)\subo$ on the odd part $F(4)\subuno$
is the tensor product of the $8$-dimensional spin representation of $B_3$
with the natural $2$-dimensional representation of $A_1$.

\item There is a unique $31$-dimensional simple classical Lie
superalgebra $G(3)$ for which $G(3)\subo$ is a
Lie algebra of type $G_2\oplus A_1$, and its
representation on $G(3)\subuno$ is the tensor product of
the $7$-dimensional irreducible representation of $G_2$
with the natural $2$-dimensional representation of $A_1$.

\item There is a one-parameter family of $17$-dimensional simple
classical Lie superalgebras $D(2,1;t)$, $t\in
\bF\setminus\{0,-1\}$, consisting of all simple Lie superalgebras
for which the even part is a Lie algebra of type $A_1\oplus
A_1\oplus A_1$, and its representation on
the odd part is the tensor product of the natural $2$-dimensional
representations of the three $A_1$-direct summands.

\end{enumerate}

It turns out that all these exceptional simple Lie superalgebras can be constructed by means of the Tits construction reviewed in Section \ref{se:Tits}, if the Jordan algebra $J$ there is replaced by suitable Jordan superalgebras.

\medskip

First, a general fact about superalgebras. Consider the Grassmann (or exterior) algebra $G$ on a countable number of variables. That is, $G$ is the unital associative algebra over $\bF$ generated by elements $e_i$, $i=1,2,\ldots$, subject to the relations $e_i^2=0$, $e_ie_j=-e_je_i$ for any $i,j$. This Grassmann algebra is naturally graded over $\bZ_2$, with all the generators $e_i$ in $G\subuno$. In this way, $G$ becomes an associative superalgebra.

Given any other superalgebra $A=A\subo\oplus A\subuno$ over $\bF$, that is, any $\bZ_2$-graded algebra, the \emph{Grassmann envelope} of $A$ is defined to be the $\bZ_2$-graded algebra
\[
G(A)=(A\subo\otimes G\subo)\oplus(A\subuno\otimes G\subuno).
\]
This is a subalgebra of the algebra $A\otimes G$ with its natural multiplication, in particular $G(A)$ is an algebra over the commutative associative ring $G\subo$.

Given any variety $\cV$ of algebras, the superalgebra $A$ is said to be a superalgebra in the variety $\cV$ if its Grassmann envelope (as an algebra over $G\subo$) is an algebra in the variety $\cV$. In particular, a superalgebra $J=J\subo\oplus J\subuno$ is a \emph{Jordan superalgebra} if $G(J)$ is a Jordan algebra over $G\subo$. (This is equivalent to the superalgebra $J$ being supercommutative and satisfying the superized versions of the identity $x^2(yx)=(x^2y)x$ and its linearizations.)

\bigskip

Let $J=J\subo\oplus J\subuno$ be now a unital Jordan superalgebra
with a normalized trace $t_J:J\rightarrow \bF$. That is, $t_J$ is a linear
map such that $t_J(1)=1$, and $t_J(J\subuno)=0=t_J\bigl((J,J,J)\bigr)$
(see \cite[\S 1]{BE03}). Then $J=\bF 1\oplus J^0$, where $J^0=\{x\in
J: t_J(x)=0\}$, which contains $J\subuno$. For $x,y\in J^0$,
$xy=t_J(xy)1+x*y$, where $x*y=xy-t_J(xy)1$ is a supercommutative
multiplication on $J^0$. Since $(J,J,J)=[l_J,l_J](J)$ is contained
in $J^0$, the subspace $J^0$ is invariant under $\inder J=[l_J,l_J]$,
the Lie superalgebra of inner derivations. (As for Jordan algebras, $l_x:y\mapsto xy$ denotes the multiplication by the element $x$, but now the Lie bracket is the bracket of the Lie superalgebra of endomorphisms of $J$, so $d_{x,y}=[l_x,l_y]=l_xl_y-(-1)^{xy}l_yl_x$ for homogeneous elements $x,y\in J$, where $(-1)^{xy}$ is $1$ if either $x$ or $y$ is even, and it is $-1$ if both $x$ and $y$ are odd.)

\medskip

Given a composition algebra $C$ and a unital Jordan
superalgebra with a normalized trace $J$, consider the superspace
\[
\calT(C,J)=\inder C\oplus(C^0\otimes J^0)\oplus\inder J,
\]
with the superanticommutative product given by exactly the same formulas in \eqref{eq:TCJproduct}.

If the Grassmann envelope $G(J)$ satisfies the Cayley-Hamilton
equation $ch_3(x)=0$ of $3\times 3$-matrices (this is the condition that reflects the fact of being of degree $3$), where
\[
ch_3(x)=x^3-3t_J(x)x^2+\Bigl(\frac{9}{2}t_J(x)^2-\frac{3}{2}t_J(x^2)\Bigr)x-
 \Bigl(t_J(x^3)-\frac{9}{2}t_J(x^2)t_J(x)+\frac{9}{2}t_J(x)^3\Bigr)1,
\]
and where we use the same notation $t_J$ to denote the natural extension of the normalized trace in $J$ to a linear form on $G(J)$ over $G\subo$,
then with the same arguments given by Tits, $\calT(C,J)$ is shown to be a Lie superalgebra (see \cite[Sections 3 and 4]{BE03}).

However, the arguments in the previous section concerning the first three rows in the Magic Square are valid too for superalgebras and show that if the dimension of the composition algebra is $1$, $2$ or $4$ (that is, if the composition algebra $C$ is associative), then $\calT(C,J)$ is a Lie superalgebra for any Jordan superalgebra $J$.

\medskip

Only a few simple Jordan superalgebras satisfy the condition on the Cayley-Hamilton equation of degree $3$. The finite dimensional simple Jordan superalgebras have been classified in \cite{KacJordan} and \cite{MZ} (see also \cite{HK} and \cite{RZ}).

\smallskip

Given a vector superspace $V=V\subo\oplus V\subuno$ endowed with a nondegenerate supersymmetric bilinear form $b$ (so that $b$ is symmetric on $V\subo$, skew-symmetric on $V\subuno$ and $b(V\subo,V\subuno)=0=b(V\subuno,V\subo)$), the Jordan superalgebra $J(V,b)$ of this form is defined as the unital superalgebra $J=J\subo\oplus J\subuno$ with $J\subo=\bF1 \oplus V\subo$, $J\subuno=V\subuno$, with the multiplication determined by $uv=b(u,v)1$ for any $u,v\in V$. This is always a simple Jordan superalgebra. Denote by $V^{0\vert 2}$ the vector superspace with $V\subo=0$ and $\dim V\subuno=2$, endowed with a nondegenerate supersymmetric bilinear form $b$ (unique up to scaling), and by $J^{0\vert 2}$ the corresponding three-dimensional simple Jordan superalgebra $J(V^{0\vert 2},b)$. This superalgebra is trivially endowed with a normalized trace form given by $t_J(1)=1$, $t_J(u)=0$ for any $u\in V^{0\vert 2}$, and it is easy to check that its Grassmann envelope satisfies the Cayley-Hamilton equation of degree $3$.

\smallskip

Among the simple Jordan superalgebras, there is a one-parameter
family of 4-dimensional algebras $D_t \ (t \neq 0)$, having
$(D_t)\subo = \bF e \oplus \bF f$ and $(D_t)\subuno = \bF x \oplus \bF y$ where
\[
\begin{gathered}
e^2 = e, \qquad f^{2} = f, \qquad ef = 0 \\
xy = e+t f \,(= -yx), \qquad ex = \frac{1}{2} x = fx, \qquad
ey = \frac{1}{2} y = fy.
\end{gathered}
\]
For $t\ne 0,-1$, $D_t$ admits a normalized trace given by:
\[
t_J(e)=\frac{t}{1+t},\quad t_J(f)=\frac{1}{1+t},\quad t_J\bigl((D_t)\subuno\bigr)=0,
\]
and it is shown in \cite[Lemma 2.3]{BE03} that the Grassmann envelope satisfies the Cayley-Hamilton equation of degree $3$ if and only if $t$ is either $2$ or $\frac{1}{2}$. (Note that $D_t$ is isomorphic to $D_{\frac{1}{t}}$.)

\smallskip

Another simple Jordan superalgebra must enter into our discussion: the $10$-dimensional Kac's superalgebra $K_{10}$. An easy way to describe this Jordan superalgebra appeared in \cite{BE02} in terms of the smaller
Kaplansky superalgebra. The tiny Kaplansky superalgebra is
the three dimensional Jordan superalgebra $K_3=K\subo\oplus K\subuno$,
with $K\subo=\bF e$ and $K\subuno =U$, a two dimensional vector space
endowed with a nonzero alternating bilinear form $(.\vert .)$, and
multiplication given by
\[
e^2=e,\quad ex=xe=\frac{1}{2}x,\quad xy=(x\vert y)e,
\]
for any $x,y\in U$. The bilinear form $(.\vert .)$ can be extended
to a supersymmetric bilinear form by means of $(e\vert
e)=\frac{1}{2}$ and $(K\subo\vert K\subuno)=0$.

The Kac Jordan superalgebra is
\begin{equation}\label{eq:K10}
K_{10}=\bF 1\oplus (K_3\otimes K_3),
\end{equation}
with unit element $1$ and product determined \cite[(2.1)]{BE02} by
\[
(a\otimes b)(c\otimes d)=(-1)^{bc}\Bigl(ac\otimes
bd-\frac{3}{4}(a\vert c)(b\vert d)1\Bigr),
\]
for homogeneous elements $a,b,c,d\in K_3$. This superalgebra is simple if the characteristic is $\ne 2,3$. In characteristic $3$, it contains the simple ideal $K_9=K_3\otimes K_3$. Besides, the Jordan superalgebra $K_{10}$  is endowed with a unique normalized trace given by $t_J(1)=1$ and $t_J(K_3\otimes K_3)=0$.

\begin{proposition}\label{pr:degree3}
Assume that the characteristic of our ground field $\bF$ is $\ne 2,3$. Then the only finite-dimensional simple
unital Jordan superalgebras $J$ with $J\subuno \neq 0$ whose Grassmann envelope $G(J)$ satisfies
the Cayley-Hamilton equation ${ch}_3(x) = 0$ relative a normalized
trace on $J$ are, up to isomorphism, $J^{0\vert 2}$, $D_2$ and, only if the characteristic of $\bF$ is $5$, $K_{10}$.
\end{proposition}
\begin{proof}
This is proved in \cite[Proposition 5.1]{BE03}, where $K_{10}$ is missing (the argument has a flaw for $K_{10}$ as the subalgebra considered there to get a contradiction is not a unital subalgebra). The situation for $K_{10}$ is settled in \cite{McC}.
\end{proof}

\bigskip

Therefore, we can plug the three superalgebras in Proposition \ref{pr:degree3} in the Tits construction $\calT(C,J)$ (and in all the rows but the last one, we can plug also the Jordan superalgebras $D_t$ for arbitrary values of $t\ne 0$), to obtain the \emph{Supermagic Rectangle} in Table \ref{ta:SuperRectangle}. The entries in the last three columns are computed in \cite{BE03} and \cite{EldModular}.

{\tiny
\begin{table}[h!]
\begin{center}
\begin{tabular}{c|cccc|ccc}
 $\calT(C,J)$&\bigstrut$H_3(\bF)$&$H_3(\bF\times \bF)$&$H_3(\Mat_2(\bF))$&$H_3(C(\bF))$&$J^{0\vert 2}$& $D_t$&$K_{10}$\\
 \cline{1-8}
 \bigstrut$\bF$&$A_1$&$A_2$&$C_3$&$F_4$&$A_1$&$B(0,1)$&$B(0,1)\oplus B(0,1)$\\
 \bigstrut$\bF\times \bF$&$A_2$&$A_2\oplus A_2$&$A_5$&$E_6$&$B(0,1)$&$A(1,0)$&$C(3)$\\
 \bigstrut$\Mat_2(\bF)$&$C_3$&$A_5$&$D_6$&$E_7$&$B(1,1)$&$D(2,1;t)$&$F(4)$\\
 \bigstrut$C(\bF)$&$F_4$&$E_6$&$E_7$&$E_8$&$G(3)$&$F(4)$ &$\calT(C(\bF),K_{10})$\\[-2pt]
 &&&&&&($t=2$)&($\charac 5$)
\end{tabular}
\bigskip\null
\end{center}
\caption{A Supermagic Rectangle}\label{ta:SuperRectangle}
\end{table}
}

\medskip

Therefore, the exceptional simple classical Lie superalgebras $D(2,1;t)$, $G(3)$ and $F(4)$ can be obtained too by means of the Tits construction.

And in characteristic $5$ there appears another Lie superalgebra, $\calT\bigl(C(\bF),K_{10}\bigr)$. This superalgebra turns out to be a new simple modular Lie superalgebra specific to this characteristic. Its dimension is $87$, its even part is the classical simple Lie algebra $B_5$, and its odd part is the spin module for the even part. This Lie superalgebra appeared for the first time in \cite{EldModular}, where all these features are proved.

In the notation of \cite{BGLgordo}, this superalgebra appears as $\mathfrak{el}(5;5)$.

\bigskip

\section{A symmetric construction of the Magic Square}\label{se:SymCons}

Given any composition algebra $C$ over $\bF$ with norm $n$ and standard
involution $x\mapsto \bar x=n(1,x)1-x$, the algebra $H_3(C,*)$ of
$3\times 3$ hermitian matrices over $C$, where $(a_{ij})^*=(\bar
a_{ji})$, is a Jordan algebra under the symmetrized product
\begin{equation}\label{eq:Jproduct}
x\circ y= \frac{1}{2}\bigl( xy+yx\bigr).
\end{equation}

Then,
\[
\begin{split}
J=H_3(C,*)&=\left\{ \begin{pmatrix} \alpha_0 &\bar a_2& a_1\\
  a_2&\alpha_1&\bar a_0\\ \bar a_1&a_0&\alpha_2\end{pmatrix} :
  \alpha_0,\alpha_1,\alpha_2\in \bF,\ a_0,a_1,a_2\in C\right\}\\[6pt]
 &=\bigl(\oplus_{i=0}^2 \bF e_i\bigr)\oplus
     \bigl(\oplus_{i=0}^2\iota_i(C)\bigr),
\end{split}
\]
where
\[
\begin{aligned}
e_0&= \begin{pmatrix} 1&0&0\\ 0&0&0\\ 0&0&0\end{pmatrix}, &
 e_1&=\begin{pmatrix} 0&0&0\\ 0&1&0\\ 0&0&0\end{pmatrix}, &
 e_2&= \begin{pmatrix} 0&0&0\\ 0&0&0\\ 0&0&1\end{pmatrix}, \\
 \iota_0(a)&=2\begin{pmatrix} 0&0&0\\ 0&0&\bar a\\
 0&a&0\end{pmatrix},&
 \iota_1(a)&=2\begin{pmatrix} 0&0&a\\ 0&0&0\\
 \bar a&0&0\end{pmatrix},&
 \iota_2(a)&=2\begin{pmatrix} 0&\bar a&0\\ a&0&0\\
 0&0&0\end{pmatrix},
\end{aligned}
\]
for any $a\in C$. Identify $\bF e_0\oplus \bF e_1\oplus \bF e_2$ with $\bF^3$ by
means of $\alpha_0e_0+\alpha_1e_1+\alpha_2e_2\simeq
(\alpha_0,\alpha_1,\alpha_2)$. Then the commutative
multiplication \eqref{eq:Jproduct} becomes:
\begin{equation}\label{eq:Jniceproduct}
\left\{\begin{aligned}
 &(\alpha_0,\alpha_1,\alpha_2)\circ(\beta_1,\beta_2,\beta_3)=
    (\alpha_0\beta_0,\alpha_1\beta_1,\alpha_2\beta_2),\\
 &(\alpha_0,\alpha_1,\alpha_2)\circ \iota_i(a)
  =\frac{1}{2}(\alpha_{i+1}+\alpha_{i+2})\iota_i(a),\\
 &\iota_i(a)\circ\iota_{i+1}(b)=\iota_{i+2}(\bar a\bar b),\\
 &\iota_i(a)\circ\iota_i(b)=2n(a,b)\bigl(e_{i+1}+e_{i+2}\bigr),
\end{aligned}\right.
\end{equation}
for any $\alpha_i,\beta_i\in \bF$, $a,b\in C$, $i=0,1,2$, and where
indices are taken modulo $3$.

Note that \eqref{eq:Jniceproduct} shows that $J$ is graded over
$\bZ_2\times \bZ_2$ with:
\[
J_{(0,0)}=\bF^3,\quad J_{(1,0)}=\iota_0(C),\quad
 J_{(0,1)}=\iota_1(C),\quad J_{(1,1)}=\iota_2(C)
\]
and, therefore, $\der J$ is accordingly graded over
$\bZ_2\times\bZ_2$:
\[
(\der J)_{(i,j)}=\{ d\in\der J: d\bigl(J_{(r,s)}\bigr)\subseteq
J_{(i+r,j+s)}\ \forall r,s=0,1\}.
\]

It is not difficult to prove (see \cite[Lemma 3.4]{CunEld2}) that
$(\der J)_{(0,0)}=\{ d\in \der J : d(e_i)=0\ \forall
i=0,1,2\}$.

\smallskip

Now,  for any  $d\in(\der J)_{(0,0)}$, there are linear maps $d_i\in \End_\bF(C)$, $i=0,1,2$ such that
$d\bigl(\iota_i(a)\bigr)=\iota_i\bigl(d_i(a)\bigr)$ for any $a\in C$
and $i=0,1,2$. Now, for any $a,b\in C$ and $i=0,1,2$:
\[
\begin{split}
0&=2n(a,b)d(e_{i+1}+e_{i+2})=d\bigl(\iota_i(a)\circ\iota_i(b)\bigr)\\
 &=\iota_i\bigl(d_i(a)\bigr)\circ\iota_i(b)
    +\iota_i(a)\circ\iota_i\bigl(d_i(b)\bigr)\\
 &=2\bigl(n(d_i(a),b)+n(a,d_i(b)\bigr)(e_{i+1}+e_{i+2}),
\end{split}
\]
so $d_i$ belongs to the orthogonal Lie algebra
$\frso(C,n)$. Also, if we write $a\bullet b=\bar a\bar b$, we have:
\[
\begin{split}
\iota_i\bigl(d_i(a\bullet b)\bigr)&=
    d\bigl(\iota_i(a\bullet b)\bigr)=d\bigl(\iota_{i+1}(a)\circ
    \iota_{i+2}(b)\bigr)\\
  &=d\bigl(\iota_{i+1}(a)\bigr)\circ \iota_{i+2}(b) +
    \iota_{i+1}(a)\circ d\bigl(\iota_{i+2}(b)\bigr)\\
  &=\iota_{i+1}\bigl(d_{i+1}(a)\bigr)\circ \iota_{i+2}(b) +
    \iota_{i+1}(a)\circ \iota_{i+2}\bigl(d_{i+2}(b)\bigr)\\
  &=\iota_i\Bigl(d_{i+1}(a)\bullet b+ a\bullet d_{i+2}(b)\Bigr),
\end{split}
\]
which shows that the triple $(d_0,d_1,d_2)$ satisfies the following condition:
\[
d_i(\bar a\bar b)=\overline{d_{i+1}(a)}\bar b+\bar a\overline{d_{i+2}(b)}
\]
for any $a,b\in C$ and $i=0,1,2$ (indices modulo $3$). One can check that the condition above for $i=0$ implies the conditions for $i=1$ and $i=2$. This leads us to the \emph{triality Lie algebra} of $C$, which is the Lie subalgebra of the direct sum of three copies of the orthogonal Lie algebra $\frso(C,n)$ given by the following equation:
\begin{equation}\label{eq:triC}
\tri(C)=\{(d_0,d_1,d_2)\in \frso(C,n)^3: d_0(\bar a\bar b)=\overline{d_{1}(a)}\bar b+\bar a\overline{d_{2}(b)}\ \forall a,b\in C\}.
\end{equation}
It turns out that $(\der J)_{(0,0)}$ is then isomorphic to the triality Lie algebra $\tri(C)$, by means of the linear map:
\[
\begin{split}
\tri(C)&\longrightarrow (\der J)_{(0,0)}\\
 (d_0,d_1,d_2)&\mapsto D_{(d_0,d_1,d_2)},
\end{split}
\]
such that
\[
\left\{\begin{aligned} &D_{(d_0,d_1,d_2)}(e_i)=0,\\
   &D_{(d_0,d_1,d_2)}\bigl(\iota_i(a)\bigr)=
   \iota_i\bigl(d_i(a)\bigr)
   \end{aligned}\right.
\]
for any $i=0,1,2$ and $a\in C$.

\smallskip

Also, for any $i=0,1,2$ and $a\in C$, consider the following inner
derivation of the Jordan superalgebra $J$:
\[
D_i(a)=2\bigl[ l_{\iota_i(a)},l_{e_{i+1}}\bigr]
\]
(indices modulo $3$) where, as before, $l_x$ denotes the multiplication by $x$
in $J$. Note that the restriction of $l_{e_i}$ to
$\iota_{i+1}(C)\oplus\iota_{i+2}(C)$ is half the identity, so the
inner derivation $\bigl[ l_{\iota_i(a)},l_{e_{i}}\bigr]$ is trivial
on $\iota_{i+1}(C)\oplus\iota_{i+2}(C)$, which generates $J$. Hence
\[
\bigl[l_{\iota_i(a)},l_{e_{i}}\bigr]=0
\]
for any $i=0,1,2$ and $a\in C$. Also, $l_{e_0+e_1+e_2}$ is the
identity map, so we have $\bigl[ l_{\iota_i(a)},l_{e_0+e_{1}+e_2}\bigr]=0$,
and hence
\[
D_i(a)=2\bigl[ l_{\iota_i(a)},l_{e_{i+1}}\bigr]=
-2\bigl[l_{\iota_i(a)},l_{e_{i+2}}\bigr].
\]

A straightforward computation with \eqref{eq:Jniceproduct} gives
\[
\begin{split}
&D_i(a)(e_i)=0,\ D_i(a)(e_{i+1})=\frac{1}{2} \iota_i(a),\
  D_i(a)(e_{i+2})=-\frac{1}{2}\iota_i(a),\\
&D_i(a)\bigl(\iota_{i+1}(b)\bigr)=-\iota_{i+2}(a\bullet b),\\
&D_i(a)\bigl(\iota_{i+2}(b)\bigr)=\iota_{i+1}(b\bullet a),\\
&D_i(a)\bigl(\iota_i(b)\bigr)=2n(a,b)(-e_{i+1}+e_{i+2}),
\end{split}
\]
for any $i=0,1,2$ and any homogeneous elements $a,b\in C$.

Denote by $D_{\tri(C)}$ the linear span of the $D_{(d_0,d_1,d_2)}$'s, $(d_0,d_1,d_2)\in \tri(C)$, and by $D_i(C)$ the linear span of the $D_i(a)$'s, $a\in C$. Then $D_{\tri(C)}=(\der J)_{(0,0)}$,
$D_0(C)=(\der J)_{(1,0)}$, $D_1(C)=(\der J)_{(0,1)}$, and
$D_2(C)=(\der J)_{(1,1)}$ (see \cite[Lemma 3.11]{CunEld2}).

Therefore, the $\bZ_2\times\bZ_2$-grading of $\der J$ becomes
\[
\der J=D_{\tri(C)}\oplus\bigl(\oplus_{i=0}^2 D_i(C)\bigr).
\]
At least as vector spaces, and assuming that the characteristic of $\bF$ is $\ne 2,3$, we obtain isomorphisms:
\begin{align*}
&J=H_3(C)\simeq
\bF^3\oplus\bigl(\oplus_{i=0}^2\iota_i(C)\bigr),\\[2pt]
&J_0\simeq \bF^2\oplus \bigl(\oplus_{i=0}^2\iota_i(C)\bigr),\\[4pt]
&\der J\simeq \tri(C)\oplus\bigl(\oplus_{i=0}^2\iota_i(C)\bigr).
\end{align*}

Whence, if $C$ is a composition algebra with norm $n$, and $J$ is the simple Jordan algebra of hermitian $3\times 3$ matrices over a second composition algebra $C'$ with norm $n'$, the Lie algebra considered in the Tits construction (note that $\inder J=\der J$) can be split into pieces and then rearranged as follows, at least as a vector space:
\begin{align*}
\calT(C,J)&=\inder C\oplus (C^0\otimes J^0)\oplus\inder J\\[4pt]
 &\simeq \der C\oplus (C^0\otimes \bF^2)\oplus\bigl(\oplus_{i=0}^2C^0\otimes
 \iota_i(C')\bigr)\oplus\bigl(\tri(C')\oplus(\oplus_{i=0}^2\iota_i(C'))\bigr)\\[4pt]
 &\simeq\bigl(\der C\oplus (C^0\otimes \bF^2)\bigr)\oplus\tri(C')\oplus
 \bigl(\oplus_{i=0}^2\iota_i(C\otimes C')\bigr)\\[4pt]
 &\simeq\bigl(\tri(C)\oplus\tri(C')\bigr)\oplus\bigl(\oplus_{i=0}^2\iota_i(C\otimes
 C')\bigr)
\end{align*}
where $\iota_i(C\otimes C')$ indicates a copy of the tensor product $C\otimes C'$, and where we have used that $\der C$ identifies canonically with the subalgebra of those elements in $\tri(C)$ with equal components, and then $\tri(C)$ decomposes as the direct sum of this subalgebra and two copies of the module for $\der C$ formed by the subspace of trace zero elements in $C$ (this follows, for instance, from the arguments in \cite[Chapter III.8]{Schafer} in dimension $8$ and \cite[Corollary 3.4]{EldIbero1} in dimension $\leq 4$).

\medskip

Then (see \cite{BS1}, \cite{BS2}, \cite{LM1}, \cite{LM2} or \cite{EldIbero1}) the Lie bracket in $\calT(C,J)$ can be transferred to the vector space
\[
\frg(C,C')=\bigl(\tri(C)\oplus\tri(C')\bigr)\oplus\bigl(\oplus_{i=0}^2
\iota_i(C\otimes C')\bigr),
\]
as follows:
\begin{equation}\label{eq:gCC'}
\begin{split}
\bullet&\ \textrm{$\tri(C)\oplus\tri(C')$ is a Lie subalgebra of $\frg(C,C')$,}\\
\bullet&\ [(d_0,d_1,d_2),\iota_i(x\otimes
 x')]=\iota_i\bigl(d_i(x)\otimes x'\bigr),\\
\bullet&\
 [(d_0',d_1',d_2'),\iota_i(x\otimes
 x')]=\iota_i\bigl(x\otimes d_i'(x')\bigr),\\
\bullet&\ [\iota_i(x\otimes x'),\iota_{i+1}(y\otimes y')]=
 \iota_{i+2}\bigl((\bar x\bar y)\otimes (\bar x' \bar y')\bigr)\ \textrm{(indices modulo $3$),}\\
\bullet&\ [\iota_i(x\otimes x'),\iota_i(y\otimes y')]=
 n'(x',y')\theta^i(t_{x,y}) +
 n(x,y)\theta'^i(t'_{x',y'}),
\end{split}
\end{equation}
for any $x,y\in C$, $x',y'\in C'$, $(d_0,d_1,d_2)\in \tri(C)$, and $(d_0',d_1',d_2')\in \tri(C')$, where
\[
t_{x,y}=\bigl(n(x,.)y-n(y,.)x,\tfrac{1}{2}n(x,y)1-R_{\bar x}R_y,\tfrac{1}{2}n(x,y)1-L_{\bar x}L_y\bigr),
\]
$L$ and $R$ denote, respectively, left and right multiplications, and $\theta$ is the automorphism of $\tri(C)$ such that $\theta\bigl((d_0,d_1,d_2)\bigr)=(d_2,d_0,d_1)$. Similarly for $t_{x',y'}$ and $\theta'$.

In this way, Tits construction becomes a construction depending symmetrically on two composition algebras.

\smallskip

There is another advantage of this symmetric construction, and it is that the above definition of the Lie bracket in $\frg(C,C')$ makes sense too in characteristic $3$. One thus obtains the Magic Square in characteristic $3$ (Table \ref{ta:MS3}).

\begin{table}[h!]
\begin{center}
\begin{tabular}{rc|ccccc}
 \multicolumn{2}{c}{}&\multicolumn{4}{c}{$\dim C'$}&\null\qquad\qquad\null\\[4pt]
 &$\frg(C,C')$&$1$&$2$&$4$&$8$\\
 \cline{2-6}
 &\bigstrut$1$&$A_1$&$\tilde A_2$&$C_3$&$F_4$\\
 &\bigstrut$2$&$\tilde A_2$&$\tilde A_2\oplus \tilde A_2$&$\tilde A_5$&$\tilde E_6$\\
 \smash{\raise 6pt\hbox{$\dim C$}}&\bigstrut$4$&$C_3$&$\tilde A_5$&$D_6$&$E_7$\\
 &\bigstrut$8$&$F_4$&$\tilde E_6$&$E_7$&$E_8$
\end{tabular}
\end{center}
\bigskip\null
\caption{Magic Square in characteristic $3$}\label{ta:MS3}
\end{table}

In Table \ref{ta:MS3}, $\tilde A_2$ (respectively $\tilde A_5$) denotes the projective general Lie algebra $\frpgl_3(\bF)$ (respectively $\frpgl_6(\bF)$), which is not simple, but contains the codimension one simple ideal $\frpsl_3(\bF)$ (respectively $\frpsl_6(\bF)$). In the same vein, $\tilde E_6$ denotes a $78$-dimensional Lie algebra which contains a unique codimension one simple ideal: the simple Lie algebra of type $E_6$ in characteristic $3$ (whose dimension is $77$!).

\bigskip

\section{Composition superalgebras}\label{se:ComSuper}

A quadratic superform on a $\bZ_2$-graded vector space
$U=U\subo\oplus U\subuno$ over our field $\bF$ is a pair
$q=(q\subo,\bil)$ where $q\subo :U\subo\rightarrow \bF$ is a quadratic
form, and
 $\bil:U\times U\rightarrow \bF$ is a supersymmetric even bilinear form
such that $\bil\vert_{U\subo\times U\subo}$ is the polar form of $q\subo$:
\[
\bil(x\subo,y\subo)=q\subo(x\subo+y\subo)-q\subo(x\subo)-q\subo(y\subo)
\]
for any $x\subo,y\subo\in U\subo$.

The quadratic superform $q=(q\subo,\bil)$ is said to be
\emph{regular} if the bilinear form $\bil$  is
nondegenerate.

\smallskip

Then a unital superalgebra $C=C\subo\oplus C\subuno$, endowed with
a regular quadratic superform $q=(q\subo,\bil)$, called the
\emph{norm}, is said to be a \emph{composition superalgebra} (see
\cite{EldOkuCompoSuper}) in case
\begin{subequations}
\begin{align}
&q\subo(x\subo y\subo)=q\subo(x\subo)q\subo(y\subo),\label{eq:qcompo1}\\
&\bil(x\subo y,x\subo z)=q\subo(x\subo)\bil(y,z)=\bil(yx\subo,zx\subo),\label{eq:qcompo2}\\
&\bil(xy,zt)+(-1)^{  x   y  +
x
 z +  y   z }\bil(zy,xt)=(-1)^{
 y   z }\bil(x,z)\bil(y,t),\label{eq:qcompo3}
\end{align}
\end{subequations}
for any $x\subo,y\subo\in C\subo$ and homogeneous elements
$x,y,z,t\in C$. (As we are working in characteristic $\ne 2$, it is enough to consider equation \eqref{eq:qcompo3}.)


\smallskip

Nontrivial composition superalgebras appear only over fields of characteristic $3$ (see \cite{EldOkuCompoSuper}):

\begin{itemize}

\item Let $V$ be a $2$-dimensional vector space over $\bF$,
endowed with a nonzero alternating bilinear form $\langle .\vert
.\rangle$ (that is, $\langle v\vert v\rangle =0$ for any $v\in V$).  Consider the superspace $B(1,2)$ (see \cite{She97}) with
\[
B(1,2)\subo =\bF 1,\qquad\text{and}\qquad B(1,2)\subuno= V,
\]
endowed with the supercommutative multiplication given by
\[
1x=x1=x\qquad\text{and}\qquad uv=\langle u\vert v\rangle 1
\]
for any $x\in B(1,2)$ and $u,v\in V$, and with the quadratic
superform $q=(q\subo,\bil)$ given by:
\[
q\subo(1)=1,\quad \bil(u,v)=\langle u\vert v\rangle,
\]
for any $u,v\in V$. If the characteristic of $\bF$ is $3$, then
$B(1,2)$ is a composition superalgebra (\cite[Proposition
2.7]{EldOkuCompoSuper}).

\smallskip

\item Moreover, with $V$ as before, let $f\mapsto \bar f$ be the
associated symplectic involution on $\End_\bF(V)$ (so $\langle
f(u)\vert v\rangle =\langle u\vert\bar f(v)\rangle$ for any $u,v\in
V$ and $f\in\End_\bF(V)$). Consider the superspace $B(4,2)$ (see
\cite{She97}) with
\[
B(4,2)\subo=\End_\bF(V),\qquad\text{and}\qquad B(4,2)\subuno=V,
\]
with multiplication given by the usual one (composition of maps) in
$\End_\bF(V)$, and by
\[
\begin{split}
&v\cdot f=f(v)=\bar f\cdot v \in V,\\
&u\cdot v=\langle .\vert u\rangle v\in \End_\bF(V)
\end{split}
\]
for any $f\in\End_\bF(V)$ and $u,v\in V$, where $\langle .\vert u\rangle v$ denotes the endomorphism $w\mapsto \langle w\vert u\rangle v$; and with quadratic superform
$q=(q\subo,\bil)$ such that
\[
q\subo(f)=\det(f),\qquad\bil(u,v)=\langle u\vert v\rangle,
\]
for any $f\in \End_\bF(V)$ and $u,v\in V$. If the
characteristic is $3$, $B(4,2)$ is a composition superalgebra
(\cite[Proposition 2.7]{EldOkuCompoSuper}).

\end{itemize}

\smallskip

Given any composition superalgebra $C$ with norm $q=(q\subo,\bil)$, its
standard involution is given by
\[
x\mapsto \bar x=\bil(x,1)1-x,
\]
and its \emph{triality Lie superalgebra} $\tri(C)=\tri(C)\subo\oplus\tri(C)\subuno$ is
defined by:
\begin{multline*}
\tri(C)_{\bar\imath}=\{ (d_0,d_1,d_2)\in\frosp(C,b)^3_{\bar\imath}:\\
d_0(x\bullet y)=d_1(x)\bullet y+(-1)^{i  x }x\bullet
d_2(y)\ \forall x,y\in C\subo\cup C\subuno\},
\end{multline*}
for $\bar \imath= \bar 0,\bar 1$, where $x\bullet y=\bar x\bar y$ for any $x,y\in C$,  and $\frosp(S,b)$ denotes the
associated orthosymplectic Lie superalgebra. The bracket in
$\tri(C)$ is given componentwise.

\smallskip

Now the construction given in Section \ref{se:SymCons} of the Lie algebra $\frg(C,C')$ is valid in this setting. Therefore, given two composition superalgebras $C$ and $C'$, consider the Lie superalgebra:
\[
\frg=\frg(C,C')=\bigl(\tri(C)\oplus\tri(C')\bigr)\oplus\bigl(\oplus_{i=0}^2
\iota_i(C\otimes C')\bigr),
\]
where $\iota_i(C\otimes C')$ is just a copy of $C\otimes C'$
($i=0,1,2$),  with bracket given by the superversion of \eqref{eq:gCC'}:

\smallskip

\begin{itemize}
\item the Lie bracket in $\tri(C)\oplus\tri(C')$, which thus becomes  a Lie subsuperalgebra of $\frg$,
\smallskip

\item $[(d_0,d_1,d_2),\iota_i(x\otimes
 x')]=\iota_i\bigl(d_i(x)\otimes x'\bigr)$,
\smallskip

\item
 $[(d_0',d_1',d_2'),\iota_i(x\otimes
 x')]=(-1)^{  d_i'   x }\iota_i\bigl(x\otimes d_i'(x')\bigr)$,
\smallskip

\item $[\iota_i(x\otimes x'),\iota_{i+1}(y\otimes y')]=(-1)^{
x'   y }
 \iota_{i+2}\bigl((x\bullet y)\otimes (x'\bullet y')\bigr)$ (indices modulo
 $3$),
\smallskip

\item $[\iota_i(x\otimes x'),\iota_i(y\otimes y')]=
 (-1)^{  x   x' +  x   y'  +
   y   y' }
 \bil'(x',y')\theta^i(t_{x,y})$ \newline \null\hspace{2.5 in} $+
 (-1)^{  y   x' }
 \bil(x,y)\theta'^i(t'_{x',y'})$,

\end{itemize}
\smallskip

\noindent
for any $i=0,1,2$ and homogeneous $x,y\in C$, $x',y'\in C'$,
$(d_0,d_1,d_2)\in\tri(C)$, and $(d_0',d_1',d_2')\in\tri(C')$. Here
$\theta$ denotes the natural automorphism
$\theta:(d_0,d_1,d_2)\mapsto (d_2,d_0,d_1)$ in $\tri(C)$, while $t_{x,y}$ is defined now by
\[
t_{x,y}=\bigl(\sigma_{x,y},\tfrac{1}{2}\bil(x,y)1-r_xl_y,\tfrac{1}{2}\bil(x,y)1-l_xr_y\bigr)
\]
with $l_x(y)=x\bullet y$, $r_x(y)=(-1)^{xy}y\bullet x$, and
\[
\sigma_{x,y}(z)=(-1)^{yz}\bil(x,z)y-(-1)^{x(y+z)}\bil(y,z)x
\]
for homogeneous $x,y,z\in S$. Also $\theta'$
and $t'_{x',y'}$ denote the analogous elements for $\tri(C')$.

\smallskip

Assuming the characteristic of $\bF$ is $3$, the Lie superalgebras $\frg(C,C')$, where $C,C'$ run over $\{\bF,\bF\times\bF,\Mat_2(\bF),C(\bF), B(1,2),B(4,2)\}$, appear in Table \ref{ta:supermagicsquare}, which has been obtained in \cite{CunEld1}, and which extends Table \ref{ta:MS3} with the addition of two rows and columns filled with Lie superalgebras.

{\small
\begin{table}[h!]
\begin{center}
\begin{tabular}{c|cccc|cc}
 $\frg(C,C')$&$\bF$&$\bF\times \bF$&$\Mat_2(\bF)$&$C(\bF)$&$B(1,2)$&$B(4,2)$\\
 \hline
 \bigstrut$\bF$&$A_1$&$\tilde A_2$&$C_3$&$F_4$&$\frpsl_{2,2}$&$\frsp_6\oplus (14)$\\
 \bigstrut$\bF\times \bF$&&$\tilde A_2\oplus \tilde A_2$&$\tilde A_5$&$\tilde E_6$&$\bigl(\frpgl_3\oplus\frsl_2\bigr)\oplus\bigl(\frpsl_3\otimes (2)\bigr)$&
    $\frpgl_6\oplus (20)$\\
 \bigstrut$\Mat_2(\bF)$&&&$D_6$&$E_7$&
 $\bigl(\frsp_6\oplus\frsl_2\bigr)\oplus\bigl((13)\otimes (2)\bigr)$
    &$\frso_{12}\oplus spin_{12}$\\
 \bigstrut$C(\bF)$&&&&$E_8$&
 $\bigl(\frf_4\oplus\frsl_2\bigr)\oplus\bigl((25)\otimes (2)\bigr)$&
      $\fre_7\oplus (56)$\\
 \hline
 \bigstrut$B(1,2)$&&&&&
 $\frso_7\oplus 2spin_7$ &$\frsp_8\oplus(40)$\\
 \bigstrut$B(4,2)$&&&&&&$\frso_{13}\oplus spin_{13}$\\
\end{tabular}
\smallskip\null
\end{center}
\caption{A Supermagic Square in characteristic
$3$}\label{ta:supermagicsquare}
\end{table}
}

Since the construction of $\frg(C,C')$ is symmetric, only the
entries above the diagonal are needed. In Table
\ref{ta:supermagicsquare}, the even and odd parts of the nontrivial superalgebras in the
table which have no counterpart in the classification in
characteristic $0$ (\cite{Kac77}) are displayed, $spin$ denotes the
spin module for the corresponding orthogonal Lie algebra, while
$(n)$ denotes a module of dimension $n$, whose precise description is given in \cite{CunEld1}. Thus, for example,
$\frg(\Mat_2(\bF)),B(1,2))$ is a Lie superalgebra whose even part is
(isomorphic to) the direct sum of the symplectic Lie algebra
$\frsp_6$ and of $\frsl_2$, while its odd part is the tensor
product of a $13$-dimensional module for $\frsp_6$ and the
natural $2$-dimensional module for $\frsl_2$.

A precise description of these modules and of the Lie superalgebras
as Lie superalgebras with a Cartan matrix is given in \cite{CunEld1}. All the
inequivalent Cartan matrices for these Lie superalgebras are listed in \cite{BGL}.

Denote by $\frg(n,m)$ the Lie superalgebra in Table \ref{ta:supermagicsquare} with $\dim C=n$ and $\dim C'=m$. Then all the Lie superalgebras $\frg(n,m)$ are simple unless either $n$ or $m$ is equal to $2$. The Lie superalgebra $\frg(2,3)$ has a unique simple ideal of codimension $1$ whose even part is $\frpsl_3\oplus\frsl_2$, while $\frg(2,6)$ has a unique simple ideal of codimension one whose even part is $\frpsl_6$.

\smallskip

In this way, a family of $10$ new simple Lie superalgebras over algebraically closed fields of characteristic $3$ appear in Table \ref{ta:supermagicsquare}.
Actually, with the exception of $\frg(3,6)$, all these superalgebras had appeared for the first time in \cite{EldNew3}, where the $\frg(n,3)$ and $\frg(n,6)$ for $n=1,2,4$, or $8$ were described in terms of some instances of the so called symplectic and orthogonal triple systems, and in \cite{EldModular}, where $\frg(3,3)$ and $\frg(6,6)$ appeared in a different way.

\bigskip


\section{Some further results and remarks}

In this section we mention several recent results related with the Tits construction and its extensions.

\subsection{Bouarroudj, Grozman and Leites classification}

The finite dimensional modular Lie superalgebras over algebraically closed fields with indecomposable symmetrizable Cartan matrices (or contragredient Lie superalgebras) have been classified in \cite{BGLgordo} under some technical hypotheses, even in characteristic $2$, where a suitable definition of Lie superalgebra is given first.

In characteristic $p\geq 3$, apart from the Lie superalgebras obtained as the analogues of the Lie superalgebras in Kac's classification in characteristic $0$ by reducing the Cartan matrices modulo $p$, there are only the following exceptions:

\begin{enumerate}

\item Two exceptions in characteristic $5$: $\frbr(2;5)$ of dimension $10\vert 12$ (that is the even part has dimension $10$ and the odd part dimension $12$), and $\frel(5;5)$ of dimension $55\vert 32$.

\item The family of exceptions given by the Lie superalgebras in the Supermagic Square in characteristic $3$ (Table \ref{ta:supermagicsquare}).

\item Another two exceptions in characteristic $3$, similar to the ones in characteristic $5$: $\frbr(2;3)$ of dimension $10\vert 8$, and $\frel(5;3)$ of dimension $39\vert 32$.
\end{enumerate}

\smallskip
Besides the superalgebras in the Supermagic Square in characteristic $3$ (Table \ref{ta:supermagicsquare}), it turns out that the superalgebra $\frel(5;5)$ is the Lie superalgebra $\calT(C(\bF),K_{10})$ in the Supermagic Rectangle (Table \ref{ta:SuperRectangle}), while the superalgebra $\frel(5;3)$ lives (as a natural maximal subalgebra) in the Lie superalgebra $\frg(3,8)$ of the Supermagic Square in characteristic $3$ (see \cite{EldModels}).

Therefore only two of the exceptions: $\frbr(2;3)$ and $\frbr(2;5)$ do not seem to be connected to an extension of the Tits construction.

The simple Lie superalgebra $\frbr(2;3)$ appeared for the first time in \cite{EldNew3}, related to an eight dimensional \emph{symplectic triple system}. Finally, a simple model of the Lie superalgebra $\frbr(2;5)$ appears in \cite{EldModels}.

\medskip

\subsection{Jordan algebras and superalgebras}

The Lie superalgebras that appear in the Supermagic Square in characteristic $3$ in Table \ref{ta:supermagicsquare} are strongly related to some Jordan algebras and superalgebras.

Thus the even part of the Lie superalgebras $\frg\bigl(B(1,2),C\bigr)$, for a composition algebra $C$, is the direct sum of the three-dimensional simple Lie algebra $\frsl_2$ and the Lie algebra of derivations of the simple Jordan algebra $J=H_3(C,*)$ of $3\times 3$ hermitian matrices over $C$, while its odd part is the tensor product of the two-dimensional natural irreducible module for $\frsl_2$ with the module for $\der J$ consisting of the quotient of the subspace $J^0$ of trace zero elements in $J$ modulo the subspace spanned by the identity matrix. (Since the characteristic is $3$, the trace of the identity matrix is $0$.)

On the other hand, the Lie superalgebras $\frg(\bigl(B(4,2),C\bigr)$, with $C$ as above, are related to the simple Freudenthal triple system defined on the subspace of $2\times 2$ matrices with diagonal entries in $\bF$ and off-diagonal matrices in $J=H_3(C,*)$. (See \cite{CunEld2} for the details.)

\smallskip

From a different perspective (see \cite{CunEldIvan}), the Lie superalgebras $\frg\bigl(B(1,2),C\bigr)$ (respectively $\frg\bigl(B(4,2),C\bigr)$), for an associative composition algebra $C$, are related to the simple Jordan superalgebra of the hermitian $3\times 3$ matrices over the composition superalgebra $B(1,2)$ (respectively $B(4,2)$). In particular, if $C=\Mat_2(\bF)$ then we get the Tits-Kantor-Koecher Lie superalgebras of these simple Jordan superalgebras.

Moreover, the simple Lie superalgebra $\frg\bigl(B(1,2),B(1,2)\bigr)$ is the Tits-Kantor-Koecher Lie superalgebra of the $9$-dimensional simple Kac Jordan superalgebra $K_9$ (recall that the $10$-dimensional Kac superalgebra in \eqref{eq:K10} is no longer simple in characteristic $3$, but contains a $9$-dimensional simple ideal, which is the tensor product (in the graded sense) of two copies of the tiny Kaplansky superalgebra $K_3$.)

\medskip

\subsection{The fourth row of Tits construction in characteristic $3$}

If the characteristic is $3$, the Lie algebra of derivations of the Cayley algebra $C(\bF)$ is no longer simple, but contains the simple ideal consisting of its inner derivations, which consists just of the adjoint maps $\ad_x:y\mapsto [x,y]$. This ideal is isomorphic to the projective special linear algebra $\frpsl_3(\bF)$ (see \cite{AMEN}). It makes sense then to consider a modified Tits construction:
\[
\begin{split}
\tilde\calT\bigl(C(\bF),J\bigr)&=\ad_{C^0}\oplus (C^0\otimes J^0)\oplus\inder J\\
&\simeq (C^0\otimes J)\oplus \inder J,
\end{split}
\]
with a bracket like the one in \eqref{eq:TCJproduct}. Then (see \cite{EldTits3}) $\tilde\calT(C,J)$ becomes a Lie algebra if and only if $J$ is a commutative alternative algebra (these two conditions imply that $J$ is a Jordan algebra).

The simple commutative alternative algebras are just the fields, so nothing interesting appears here.

But there are two types of simple commutative alternative superalgebras
\cite{She97} in characteristic $3$, apart from the fields. The easiest simple commutative alternative superalgebra is the Jordan superalgebra $J^{0\vert 2}$ which has already appeared in the construction of the Supermagic Rectangle (Table \ref{ta:SuperRectangle}) in Section \ref{se:JordanSuperalgebras}.

It follows that the simple Lie superalgebra $\tilde\calT(C(\bF),J^{0\vert 2})$ coincides with the unique simple ideal of codimension one in the Lie superalgebra $\frg(2,3)=\frg(\bF\times\bF,B(1,2))$ in the Supermagic Square in characteristic $3$ (Table \ref{ta:supermagicsquare}), so nothing new appears here.

However, the other family of simple commutative alternative superalgebras in characteristic $3$ consists of the superalgebras $B=B(\Gamma,d)$, where:

\begin{itemize}
\item $\Gamma$ is a commutative associative algebra,

\item $d\in\der\Gamma$ is a derivation such that there is no proper ideal of $\Gamma$ invariant under the action of $d$ (that is, $\Gamma$ is  $d$-simple),

\item $a(bu)=(ab)u=(au)b$,\quad $(au)(bu)=ad(b)-d(a)b$, for any
$a,b\in \Gamma$.
\end{itemize}

\smallskip

The most natural example of this situation is obtained if $\Gamma$ is taken to be the algebra of divided powers
\[
\calO(1;n)=\espan{t^{(r)}: 0\leq r\leq 3^n-1},
\]
where $t^{(r)}t^{(s)}=\binom{r+s}{r}t^{(r+s)}$, and with $d$ the derivation given by $d(t^{(r)})=t^{(r-1)}$.

The Lie superalgebras $\tilde\calT\bigl(C(\bF),B(\calO(1;n),d)\bigr)$ ($n\geq 1$) are then  simple Lie superalgebras with no counterparts in Kac's classification. These simple Lie superalgebras have appeared first in \cite{BGL} in a completely different way, and were denoted by $\textrm{Bj}(1;n\vert 7)$. (See \cite{EldTits3} for the details.)

\medskip

\subsection{Symmetric composition (super)algebras}

In Section \ref{se:SymCons} a symmetric construction of the Magic Square have been reviewed. It depends on two composition algebras, and uses their triality Lie algebras. However, nicer formulas for triality are obtained in one takes the so called \emph{symmetric composition algebras} instead of the traditional unital composition algebras (see \cite[Chapter VIII]{KMRT}).

An algebra $S$, with product denoted by $*$, endowed with a regular quadratic form $n$ permitting composition: $n(x*y)=n(x)n(y)$, and satisfying that $n(x*y,z)=n(x,y*z)$ for any $x,y,z\in S$, is called a \emph{symmetric composition algebra}.

The easiest way to obtain a symmetric composition algebra is to start with a classical composition algebra $C$ and consider the new multiplication defined by $x\bullet y=\bar x\bar y$. Then $C$, with its norm but with this new product, is a symmetric composition algebra, called a \emph{para-Hurwitz algebra}.

Actually, if one looks at the definition of the triality Lie algebra $\tri(C)$, one may notice that, in terms of this new multiplication, the definition in  \eqref{eq:triC} becomes:
\[
\tri(C)=\{(d_0,d_1,d_2)\in\frso(C,n): d_0(x\bullet y)=d_1(x)\bullet y+x\bullet d_2(y)\ \forall x,y\in C\}.
\]
(This has already been used in the definition of the triality Lie superalgebra in the previous section.) Also, the Lie bracket in $\frg(C,C')$ is better expressed in terms of the para-Hurwitz product.

The interesting fact about symmetric composition algebras is that, apart from the para-Hurwitz algebras, there appear the so called \emph{Okubo algebras}, which were introduced, under a different name, in \cite{Oku78}.

Over our algebraically closed ground field $\bF$ there is a unique Okubo algebra. If the characteristic of $\bF$ is $\ne 2,3$, this Okubo algebra is defined as the subspace $\frsl_3(\bF)$ of trace zero $3\times 3$ matrices, with the multiplication given by:
\[
x*y=\mu xy+(1-\mu)yx-\frac{1}{3}\tr(xy)1,
\]
where $\mu$ is a solution of the equation $3X(1-X)=1$, and the norm is given by $n(x)=\frac{1}{6}\tr(x^2)$. In characteristic $3$ a different definition is required (see \cite{EP96} and the references there in).

\smallskip

Our definition of the Lie algebra $\frg(C,C')$ is valid for $C$ and $C'$ being symmetric composition algebras. For para-Hurwitz algebras this is the same construction, but Okubo algebras introduce some extra features. Actually, over algebraically closed fields, no new simple Lie algebras appear, but some new properties of the simple exceptional Lie algebras can be explained in these new terms. For instance, some interesting gradings on these algebras appear naturally induced from natural gradings on Okubo algebras (see \cite{EldGradSym}).

\bigskip
\bigskip

In conclusion, I hope to have convinced the reader that the beautiful construction given by Tits \cite{Tits66} of the exceptional simple Lie algebras, has continued to be, more than forty years after its publication, a source of inspiration for further work.

\bigskip


\def\cprime{$'$}
\providecommand{\bysame}{\leavevmode\hbox
to3em{\hrulefill}\thinspace}

\end{document}